\documentclass[11pt]{article}
\usepackage{latexsym}
\usepackage{amscd}
\usepackage{amssymb}
\usepackage{amsmath}
\usepackage{amsthm}
\usepackage{enumerate}
\usepackage{verbatim}
\def\co{\colon\thinspace}
\def\d{\delta}
\def\e{\epsilon}

\newtheorem{thm}{Theorem}[section]
\newtheorem{cor}[thm]{Corollary}

\newtheorem{lem}[thm]{Lemma}
\newtheorem{prop}[thm]{Proposition}
\newtheorem{defn}[thm]{Definition}

\newtheorem{Example}[thm]{Example}
\newenvironment{ex}{\begin{Example}\rm}{\end{Example}}
\newtheorem{remark}[thm]{Remark}
\newenvironment{rmk}{\begin{remark}\rm}{\end{remark}}
\newtheorem{Fact}[thm]{Fact}

\newtheorem{Main Observation}[thm]{Main Observation}
\newenvironment{main observation}
{\begin{Main Observation}\rm}{\end{Main Observation}}
\newtheorem{Convention}[thm]{Convention}

\begin{document}
\abovedisplayskip=6pt plus3pt minus3pt
\belowdisplayskip=6pt plus3pt minus3pt
\title{\bf Finiteness theorems for nonnegatively curved vector bundles\rm}
\author{Igor Belegradek \and Vitali Kapovitch}
\date{}
\maketitle
\begin{abstract}
We prove several finiteness theorems for the
normal bundles to souls in nonnegatively curved manifolds.
More generally, we obtain finiteness results for open
Riemannian manifolds whose topology is concentrated
on compact domains of ``bounded geometry''.

\end{abstract}
\section{Introduction}
\label{Introduction}
Much of the recent work in Riemannian geometry was centered around
finiteness and precompactness theorems for various classes of
Riemannian manifolds. Some versions of precompactness results
typically work for compact domains in Riemannian manifolds. The
main point of the present paper is that one can sometimes get
diffeomorphism finiteness for ambient Riemannian manifolds
provided their topology is concentrated on a compact domains of
``bounded geometry''. We postpone the discussion of our main
technical results till section~\ref{Main observation} and
concentrate on applications to nonnegative curvature.

Recall that according to the soul theorem of J.~Cheeger and
D.~Gromoll a complete open manifold of nonnegative sectional
curvature is diffeomorphic to the total space of the normal bundle
of a compact totally geodesic submanifold which is called the
soul. One of the harder questions in the subject  is what kind of
normal bundles can occur. See~\cite{Che,Ri,Yan,GZ, GZ2} for
examples of open nonnegatively curved manifolds, and~\cite{OW,
BK2, BK3} for known obstructions. Here is our first result.

\begin{thm}\label{intro: fixed soul}
Given a closed Riemannian manifold
$S$ with $\sec(S)\ge 0$, and positive $D$, $r$, $v$, $n$,
there exists a finite collection of vector bundles over $S$
such that, for every complete open Riemannian $n$-manifolds $N$
with $\sec(N)\ge 0$ and an isometric embedding
$e\co S\to N$ of $S$ onto a soul of $N$ the normal bundle $\nu_e$
is isomorphic to a bundle of the collection provided
$e$ is homotopic to a map $f$ such that
$\mathrm{diam}(f(S))\le D$ and $\mathrm{vol}_{N}(B(p, r))\ge v$
for some $p\in f(S)$.
\end{thm}

There is also a ``variable base version'' of~\ref{intro: fixed soul}.
We say that two vector bundles $\xi$, $\xi^\prime$ over
different bases $B$, $B^\prime$ are {\it topologically equivalent}
if there is a homeomorphism  $h\co B^\prime\to B$ such that
$h^\#\xi$ is isomorphic to $\xi^\prime$. Moreover, if $h$ is a
diffeomorphism, we say that $\xi$ and $\xi^\prime$ are {\it
smoothly equivalent}.
For bundles over manifolds of dimension $\le 3$,
topological equivalence implies smooth equivalence because in
this case any homeomorphism is homotopic to a (nearby)
diffeomorphism~\cite{Mun, Moi}. The diffeomorphism type of
the total space of a vector bundle (of positive rank over a closed
manifold) is determined up to finite ambiguity by the topological
equivalence class of the bundle (see~\cite{HM, KS} if
dimension of the total space is $\ge 5$ and~\cite{Mun, Moi}
otherwise).
In the appendix we discuss to what extent
the total space determines a vector bundle and give an example of infinitely
many pairwise topologically nonequivalent vector bundles over a
closed manifold with diffeomorphic total spaces.

The following result  can be thought of as generalizations of
the Grove-Petersen-Wu finiteness theorem.

\begin{thm} \label{intro: variable soul}
Given positive $D$, $r$, $v$, $v^\prime$, $n$,
there exists a finite collection of vector bundles such that,
for any open complete Riemannian $n$-manifold $N$ with $\sec(N)\ge 0$
and a soul $S\subset N$, the normal bundle to the soul is
topologically equivalent to a bundle of the collection provided
$\mathrm{vol}(S)\ge v^\prime$ and
the inclusion $S\hookrightarrow N$
is homotopic to a map $f$ such that
$\mathrm{diam}(f(S))\le D$ and $\mathrm{vol}_{N}(B(p, r))\ge v$
for some $p\in f(S)$.
\end{thm}

We suspect that in~\ref{intro: variable soul} the lower volume
bound on the soul follows from the rest of the assumptions, and
thus can be omitted. For example, it follows from~\cite{HK} that
lower volume bound for a soul $S$ comes from a lower volume bound
on an ambient manifold $N$, that is $\mathrm{vol}_{N}B(p,r)\ge v$
implies $\mathrm{vol}(S)\ge v^\prime$ provided the distance from
$p$ to $S$ is uniformly bounded. (The latter can be also forced by
purely topological assumptions on $S$ as we show in~\ref{appl:
euler class and no collapse on soul}.) Thus we deduce the following.

\begin{cor} \label{intro: no collapse on soul}
Given positive $D$, $r$, $v$, $n$, there exists a finite
collection of vector bundles such that, for any open complete
Riemannian $n$-manifold $N$ with $\sec(N)\ge 0$ and a soul
$S\subset N$, the normal bundle to the soul is topologically
equivalent to a bundle of the collection provided
$\mathrm{diam}(S)\le D$ and $\mathrm{vol}_{N}(B(p, r))\ge v$ for
some $p\in N$ such that the distance from $S$ to $p$ is uniformly
bounded.
\end{cor}

The finiteness of homeomorphism types of total spaces
in~\ref{intro: no collapse on soul} can be easily obtained from
the parametrized version of Perelman's Stability
theorem~\cite{Per1} and the regularity properties of the distance
function; however the conclusion of~\ref{intro: no collapse on soul}
is strictly stronger (cf.~\ref{append: ex}).

There is a version of the above corollary for totally geodesic
submanifolds in Riemannian manifolds
with lower sectional curvature bound.
Another version works when ambient
manifolds have lower bound on Ricci curvature and injectivity radius.

Perelman proved in~\cite{Per2} that the distance nonincreasing
retraction onto the soul introduced in~\cite{Sha} is a
$C^{1,1}$-Riemannian submersion. We observe that a local bound on
the vertical curvatures of this submersion gives a lower volume
bound for the ambient nonnegatively curved manifold. In
particular, we deduce the following.

\begin{cor}\label{intro: fixed soul, normal curv}
Given positive $n$, $r$, $K$, and a
closed nonnegatively curved Riemannian manifold $S$,
there is a finite collection of vector bundles over $S$
such that, for any open complete Riemannian $n$-manifold
$N$ with $\mathrm{sec}(N)\ge 0$ and an isometric embedding
$e\co S\to N$ of $S$ onto a soul of $N$, the normal bundle
$\nu_e$ is isomorphic to a bundle of the collection
provided $\mathrm{diam}(S)\le D$ and there is a point
$p\in e(S)$ such that all the vertical curvatures at the points of
$B(p,r)$ are bounded above by $K$.
\end{cor}
Note that~\ref{intro: fixed soul, normal curv} generalizes the
main result of~\cite{GW} (cf.~\cite{Ta}) where the same statement
is proved for a soul isometric to the round sphere. Again, there
is a ``variable base'' version of~\ref{intro: fixed soul, normal
curv} when souls vary in  the Grove-Petersen-Wu class.
Similarly,~\ref{intro: fixed soul, normal curv} holds
when the soul varies in Cheeger-Andersen class~\cite{AC}:
$\mathrm{diam}(S)\le D$ and $\mathrm{injrad}(S)\ge i_0$
where the conclusion of finiteness up to
topological equivalence gets improved to
the finiteness up to smooth equivalence.

There is a counterpart of~\ref{intro: variable soul}
in nonpositive curvature.
Let $\sec(N)\le 0$ and $e\co M\to N$ be a
totally geodesic embedding which is a homotopy equivalence.
Then the orthogonal projection $N\to e(M)$
is distance nonincreasing. Moreover, $\mathrm{inj}(N)=\mathrm{inj}(M)$
and we obtain the following corollary (which also has
a ``fixed base'' version).

\begin{cor}\label{intro: variable nonpositive base}
Given positive $D$, $r$, $v$, $n$, $K$, there exists a finite
collection of vector bundles such that, for any totally geodesic
embedding $e\co M\to N$ of a closed Riemannian manifold $M$ into
an open complete Riemannian $n$-manifold $N$ with $\sec(N)\le 0$,
the normal bundle $\nu_e$ is topologically equivalent to a bundle
of the collection provided $\sec(M)\ge -1$, $\mathrm{vol}(M)\ge v$
and $e$ is a homotopy equivalence homotopic to a map
$f$ with $\mathrm{diam}(f(M))\le D$
such that the sectional curvature at any point of the
$r$-neighborhood of $f(M)$ is $\ge -K$.
\end{cor}

Normal bundles to totally geodesic embeddings in nonpositively
curved manifolds can be fairly arbitrary as the following
example shows.
M.~Anderson proved that
the total space $E(\xi)$ of any vector bundle $\xi$ over a closed
nonpositively curved manifold $M$ carries a complete metric with
$-1\le \sec\le 0$~\cite{And1}. Let $M$ be a closed locally
symmetric manifold of nonpositive curvature and $rank(M)\ge 2$
such that no finite cover of $M$ splits as a Riemannian product.
Let $\xi$ be an orientable vector bundle over $M$ with nonzero Euler class.
Then according to \cite[p326]{SY} the zero section $M\to E(\xi)$ is
homotopic to a harmonic map which,
by the harmonic map superrigidity~\cite{MYY}, is a totally geodesic
embedding
(after rescaling the metric on $E(\xi)$).
The normal bundle to this
totally geodesic embedding is stably isomorphic to $\xi$, and
furthermore it has the same Euler class as $\xi$.

One may wonder when there are {\it infinitely} many vector bundles
of rank $m$ over a given base $M$. For example, if $m\ge\dim(M)$,
this happens whenever $M$ has nonzero Betti number in a dimension
divisible by $4$ (e.g. if $M$ is a closed orientable manifold of
dimension divisible by $4$). The reason is that the Pontrjagin
character defines an isomorphism of $\oplus_{i>0}H^{4i}(M,\mathbb
Q)$ and $\widetilde{K}(M)\otimes\mathbb Q$ where
$\widetilde{K}(M)$ is the group of stable equivalence classes of
vector bundles over $M$. Furthermore, the Euler class defines a
one-to-one correspondence between the set of isomorphism classes
of  oriented rank $2$ bundles over $M$ and $H^2(M,\mathbb Z)$.
Also, if $M$ is a closed, orientable, and $2n$-dimensional, then
there are infinitely many rank $2n$ bundles over $M$ obtained as
pullbacks of $TS^{2n}$ via maps $M\to S^{2n}$ of nonzero degree.
%

The structure of the paper is as follows. Section~\ref{AA} reviews
some well-known results on homotopy count of maps in
equicontinuous families. Section~\ref{sec: local} discusses local
versions of precompactness theorems in~\cite{AC} and~\cite{Per1}.
The~\ref{Invariants of maps}th section provides a background in
characteristic classes and related invariants of maps.
Main technical results are proved in  section~\ref{Main observation}.
In section~\ref{sec: appl} we prove applications to
nonnegatively/nonpositively curved manifolds.
In the appendix we explain to what extent a vector bundle is determined
by its total space.

We are grateful to M.~Anderson for an illuminating communication
on the local version of~\cite{AC} and to S.~Weinberger for the
idea of~\ref{append: ex}. The first author is thankful to A.~Nicas
and I.~Hambleton for several helpful discussions on
self-equivalences of manifolds. The second author is grateful to
Kris Tapp for bringing to his attention the idea of bounding
homotopy types of maps using equicontinuity and for many helpful
conversations on nonnegatively curved manifolds.

\section{Equicontinuity and homotopy count of maps}
\label{AA}

\begin{defn} A family of maps of metric spaces 
$f_{\alpha}\co X_{\alpha}\to Y_{\alpha}$
is called $\e${\it -equicontinuous} if there exist $\d>0$ such
that $d_{Y_{\alpha}}(f_{\alpha}(x),f_{\alpha}(x^\prime))<\e$, for
any $\alpha$ and any $x,x^\prime\in X_{\alpha}$ with
$d_{X_{\alpha}}(x,x^\prime)<\d$. A family $f_{\alpha}$ is called
{\it equicontinuous} if it is $\e${\it -equicontinuous} for every $\e$. 
A family $f_{\alpha}$ is called {\it almost equicontinuous}
if for any $\e$ there exists a finite subset
$S_\e\subset\{f_{\alpha}\}$ whose complement is $\e${\it
-equicontinuous}.
\end{defn}

\begin{ex}
Assume $f_{\alpha}\co X_{\alpha}\to Y_{\alpha}$ is a family of
maps of metric spaces.

(1) If each $f_{\alpha}$ is $(\alpha, L)$-H\"older (i.e.
$d_{Y_{\alpha}}(f_{\alpha}(x),f(x^\prime))\le
Ld_{X_{\alpha}}(x,x^\prime)^\alpha$), then $\{f_{\alpha}\}$ is
equicontinuous.

(2) If each $f_{\alpha}$ is an $\e_{\alpha}$-Hausdorff
approximation (or more generally, if $f_{\alpha}$ satisfies
$$|d_{Y_{\alpha}}(f_{\alpha}(x),f_{\alpha}(x^\prime))-d_{X_{\alpha}}
(x,x^\prime)|\le\e_{\alpha}$$
for any $x,x^\prime\in X_{\alpha}$) and $\e_{\alpha}\to 0$ then,
$\{f_{\alpha}\}$ is almost equicontinuous.


(3) If $\{f_{\alpha}\}$ is almost equicontinuous and $g_{\alpha}$
is $\e_{\alpha}$-close to $f_{\alpha}$
with $\e_{\alpha}\to 0$, then $\{g_{\alpha}\}$ is almost
equicontinuous.

(4) If $\{f_{\alpha}\}$ where $f_{\alpha}\co X_{\alpha}\to
Y_{\alpha}$ is almost equicontinuous and $\{g_{\alpha}\}$ with
$g_{\alpha}\co Y_{\alpha}\to Z_{\alpha}$ is almost equicontinuous,
then $\{g_{\alpha}\circ f_{\alpha}\}$ is almost equicontinuous.
%
\end{ex}

The importance of the following result
in Riemannian geometry was first observed by M.~Gromov~\cite{GLP}.

\begin{prop}\label{AAthm}  Let $Y$ be a compact metric space such that
there exists an $\e=\e(Y)$ with the property that any two
$4\e$-close continuous maps of a compact metric space into $Y$ are
homotopic. Then, given a compact metric space $X$, any
$\e$-equicontinuous family of maps $f_{\alpha}\co X\to Y$ falls
into finitely many homotopy classes.
\end{prop}
\begin{proof}
Fix a $\d>0$ such that
$d_{Y_{\alpha}}(f_{\alpha}(x),f_{\alpha}(x^\prime))<\e$, for any
$k$ and any $x,x^\prime\in X_{\alpha}$ with
$d_{X_{\alpha}}(x,x^\prime)<\d$. Find a finite $\d/2$-net $N_X$ in
$X$ and a finite $\e$-net $N_Y$ in $Y$. Any map $f\co X\to Y$
produces a (nonunique) map $\hat f\co N_X\to N_Y$ defined so that
$\hat f(x)$ is a point of $N_Y$ whose distance to $f(x)$ is
$\le\e$. Now if $f$ and $g$ are $\e$-equicontinuous maps with
$\hat f=\hat g$, then $f$ and $g$ are $4\e$-close, hence
homotopic. In particular, $\{f_{\alpha}\}$ fall into at most
$\mathrm{card}(N_Y)^{\mathrm{card}(N_X)}$ homotopy classes.
\end{proof}
\begin{rmk}\label{rmk after AA thm}
Such an $\e(Y)$ exists if, for example, the compact metric
spaces $X$ and $Y$ are separable, finite-dimensional
ANR~\cite{PPet3}. Note that for compact, separable, finite-dimensional
metric spaces being ANR is equivalent to being locally
contractible~\cite[V.10.4]{Bor}; any such space is homotopy equivalent
to a finite cell complex~\cite{West}.
\end{rmk}

\section{Local convergence results}
\label{sec: local}

In this section we discuss local versions of the
$C^\alpha$-precompactness theorem of M.~Anderson and
J.~Cheeger~\cite{AC} and Perelman's stability theorem~\cite{Per1}.
The results provide sufficient conditions under which a sequence of
compact domains in Riemannian manifolds has uniformly bounded
geometry in the sense defined below.
In fact, the theorems in~\cite{AC} and~\cite{Per1} are
stated in a local form so we just give details needed for
``compact domains version''.

Let $U_{\alpha}$ be a family of compact domains (i.e. compact
codimension zero topological submanifolds) of Riemannian
$n$-manifolds $N_{\alpha}$. We say that $\{U_{\alpha}\}$ has {\it
uniformly bounded geometry} if any sequence of domains in the
family has a subsequence $\{U_k\}$ such that there exists a metric
space $V$, and homeomorphisms $h_k\co V\to V_k$ of $V$ onto
compact domains $V_k\supset U_k$ such that both $\{h_k\}$ and
$\{h_k^{-1}\}$ are almost equicontinuous. In case $\partial
U_k=\emptyset$, we necessarily have $U_k=V_k=N_k$.

Throughout the paper we always denote the
closed $\e$-neighborhood of a subspace $S$ by $S^\e$.
%
%
\begin{thm}\bf \cite{AC}\it\
\label{local: anderson-cheeger}
Given $\e>0$,
let $U_k$ be a sequence of compact domains with smooth boundaries
in Riemannian $n$-manifolds $N_k$
such that the closed $\e$-neighborhood $U_k^{\e}$
of $U_k$ is compact. Assume that for some positive $H$,  $V$, $i_0$,
the following holds:
$\mathrm{Ric}(U_k^\e)\ge -(n-1)H$, $\mathrm{vol}(U_k^{\e})\le V$, and
$\mathrm{inj}_{N_k}(x)\ge i_0$ for any $x\in U_k^{\e/2}$.
Then, after passing to a subsequence,
there are compact domains $V_k$ with
$U_k\subset V_k\subset U_k^{\e/2}$, a manifold $V$, and
$C^\infty$-diffeomorphisms
$h_k\co V\to V_k$ such that the pullback metrics $h_k^\ast g_k$
converge in a $C^\alpha$-topology
to a $C^\alpha$-Riemannian metric on the interior of $V$.
\end{thm}
\begin{proof} For reader's convenience we review
the argument in~\cite{AC} emphasizing its local nature. It is
proved in~\cite[pp269--270]{AC} that any domain  $U_k^{\e/2}$ as
above has an atlas of harmonic coordinate charts $F_\nu\co
B(x,r_h)\to\mathbb R^n$ where $B(x,r_h)$ is a metric ball at $x\in
U_k^{\e/2}$ whose radius $r_h\le \e/10$ depends only on the
initial data. Further, the metric tensor coefficients in the
charts $F_\nu$ are controlled in $C^\alpha$ topology. An elliptic
estimate then shows that the transition functions $F_\mu\circ
F_\nu^{-1}$ are controlled in $C^{1,\alpha}$ topology. All these
results are stated and proved locally.

Next, the relative volume comparison implies that one can choose a
finite subatlas so that there is a uniform bound on the
multiplicities of intersections of the coordinate charts and the
balls $B(x,r_h/2)$ still cover $U_k^{\e/2}$ (this argument
involves only small balls and hence is local). The lower
injectivity radius bound gives a lower bound for the volume of any
small ball that depends only on the radius of the ball~\cite{Cro}.
This, together with a upper bound on
$\mathrm{vol}(U_k^{\e})$, implies an upper bound on the number of
coordinate charts.

Finally, following Cheeger's thesis (as outlined in~\cite[pp266--267]{AC})
one can ``glue the charts together'' which proves the
theorem. Alternatively, one can follow (almost word by word) the
argument in~\cite[pp464--466]{And2} where a ``compact domain
version'' of Cheeger-Gromov convergence theorem is proved.
\end{proof}

\begin{rmk} There are many other convergence theorems,
notably those involving integral curvature bounds (see~\cite{PPet2}).
For example, Hiroshima~\cite{Hi} generalized~\cite{AC}
replacing a lower Ricci curvature bound by an integral
bound on an eigenvalue of the Ricci curvature. Hiroshima's proof
is given for complete manifolds; however, a local version of~\cite{Hi}
is likely to hold. We leave this matter for an
interested reader to clarify.
\end{rmk}

Before starting the proof of theorem~\ref{local: perelman}, we
need a local version of Packing Lemma that ensures
Gromov-Hausdorff convergence.

We say that a metric space $(X,d)$ is {\it locally interior}
if for any point $x\in X$ there exists an $\e>0$ such that
for any $y,z\in B(x,\e)$ we have $d(y,z)=\inf_\gamma L(\gamma)$
where the infimum is taken over all paths $\gamma$ connecting $y$ and $z$.
For example, all Riemannian manifolds are locally interior.
\begin{rmk}\label{int}Notice that for locally compact metric spaces the
property of being locally interior
is easily seen to be equivalent to the following one. For any
point $x\in X$ there exists an $\e>0$ such that for any $y,z\in
B(x,\e)$ there exists a sequence $p_n\in B(x,2\e)$ such that
$d(p_n,y)\to d(y,z)/2$ and $d(p_n,z)\to d(y,z)/2$.
\end{rmk}

Here is how locally interior spaces arise in this paper. Let $V_k$
be a sequence of compact domains in Riemannian $n$-manifolds.
Equip $V_k$ with induced Riemannian metrics and assume that $V_k$
converges in Gromov-Hausdorff topology to a compact metric space
$V$. Consider $f_k\co V_k\to\mathbb{R}$ defined by
$f_i(x)=dist(x,\partial V_i)$. Then each $f_i$ is 1-Lipschitz and
by Arzela-Ascoli Theorem this sequence converges to 1-Lipschitz
function $f\co V\to\mathbb{R}$. We call the open set $\{x\in U\co
f(x)>0\}$ the {\it interior} of $U$. Then it is easy to show that
the interior of $U$ is a locally interior space.

\begin{lem}\label{local: packing lemma} Given $\e>0$,
let $U_k$ be a sequence of compact connected domains with smooth
boundaries in Riemannian $n$-manifolds $N_k$ such that the closed
$\e$-neighborhood $U_k^{\e}$ of $U_k$ is compact.
Assume that for some positive $H$,  $V$, $v_0$, $r_0<\e/10$, the
following holds: $\mathrm{Ric}(U_k^\e)\ge -(n-1)H$,
$\mathrm{vol}(U_k^{\e})\le V$, and $\mathrm{vol}(B(x,r_0))\ge v_0$
for any $x\in U_k^{\e/2}$. Then, after passing to a subsequence,
the compact domains $U_k^{\e/2}$ converge in  Gromov-Hausdorff
topology to a compact metric space $U$ whose interior is a locally
interior metric space.
\end{lem}
\begin{proof} Take an arbitrary $r<r_0$.
To prove precompactness in Gromov-Hausdorff topology it is enough
to show that the number of elements in a maximal
$r$- net in $U_k^{\e/2}$ is bounded above by some
$N(r)$ independent of $k$.

Fix a maximal $r$-separated nets $N_k$ in $U_k^{\e/2}$
so that $r$-balls with centers in $N_k$ are disjoint and
$2r$-balls cover $U_k^{\e}$.
The relative volume comparison gives a uniform lower bound
for the volume of the $r$-ball centered at any point of $U_k^{\e/2}$;
say $\mathrm{vol}(B(x, r)\ge v$.
Then $\#N_k\le V/v$ and $U_k^{\e/2}$
converge in the Gromov-Hausdorff topology to
a compact metric space $U$. As we explained above
the interior of $U$ is necessarily locally interior.
\end{proof}

\begin{thm}\label{local: perelman}
\bf \cite{Per1}\it\ Given $\e>0$, let $U_k$ be a sequence of
compact connected domains with smooth boundaries in Riemannian
$n$-manifolds $N_k$ such that the closed $\e$-neighborhood of
$U_k$, denoted by $U_k^{\e}$ is compact. Assume that for some
positive $K$,  $V$, $v$, $r<\e/10$, $\sec(U_k^\e)\ge -K$,
$\mathrm{vol}(U_k^{\e})\le V$, and $\mathrm{vol}(B(x,r))\ge v$ for
any $x\in U_k^{\e/2}$. Then, after passing to a subsequence, there
are compact domains $V_k$ with $U_k^{\e/4}\subset V_k\subset
U_k^{\e/2}$, a manifold $V$, and homeomorphisms $h_k\co V\to V_k$
which are $\e_k$-Hausdorff approximations with $\e_k\to 0$ as
$k\to\infty$.
\end{thm}
\begin{proof} By~\ref{local: packing lemma} $U_k^{\e/2}$
subconverges in the Gromov-Hausdorff topology to a compact metric
space $U$ whose interior $int(U)$ is a locally interior metric
space. We are in position to apply Perelman's stability
theorem~\cite{Per1} which asserts that $int(U)$ is a topological manifold,
and moreover, any compact subset of $int(U)$ lies in a compact
domain $V\subset int(U)$ such that there are topological embedding
$h_k\co V\to U_k^{\e/2}$. Furthermore, $h_k$ induce Hausdorff
approximations which become arbitrary close to the given Hausdorff
approximations between $U$ and $U_k^{\e/2}$. Choosing $V$ large
enough, one can ensure that $h_k(V)\supset U_k^{\e/4}$ as
promised.
\end{proof}

\begin{rmk} Let $U_k$ be a sequence of compact
domains with smooth boundaries
in Riemannian $n$-manifolds $N_k$
such that each $U_k^\e$ is contained in a compact
metric ball $B(p_k, R)\subset N_k$ for some $R>\e>0$.
Assume that  $\mathrm{Ric}(B(p_k,R))\ge-(n-1)H$ for some $H>0$.
Then the absolute volume comparison implies that
$\mathrm{vol}(U_k^\e)$ is uniformly bounded
above by $B^H(R)$. Now if $\mathrm{vol}(B(x_k,\e/2))$
is uniformly bounded below, for {\it some} $x_k\in B(p_k,R-\e)$, then
the relative volume comparison ensures that
$\mathrm{vol}(B(x,r))$ is uniformly bounded below
for {\it any} $x\in B(p_k,R-\e)$ and any $r<\e/2$.

In particular,
if each $N_k$ is complete and
$\sec(N_k)\ge -1$, then any sequence of compact domains
$U_k$ has uniformly bounded geometry provided
$\mathrm{diam}(U_k)\le D$ and there are points
$x_k\in U_k$ with $\mathrm{vol}(B(x_k,r))\ge v$.
\end{rmk}

\begin{rmk}
Let $U_k$ be a sequence of compact connected domains with smooth
boundaries in Riemannian $n$-manifolds $N_k$ such that $U_k^{\e}$
is compact. Assume that $\mathrm{Ric}(U_k^\e)\ge -(n-1)H$, for
some $H>0$. Then the following two conditions are equivalent

(i) $\mathrm{vol}(U_k^{\e})\le V$, and
$\mathrm{vol}(B(x,r))\ge v$ for any $x\in U_k^{\e/2}$

(ii) there is a point $x_0\in U_k^{\e/2}$ such that
$\mathrm{vol}(B(x_0,r)\ge v$ and $\mathrm{diam}^{int} (U_k^{\e/2})\le D$
where the diameter
is taken with respect to the intrinsic distance
induced by the Riemannian metric on $U_k^{\e/2}$.

Indeed, let us show that (i)$\Rightarrow$(ii). Using  the relative
volume comparison we can make $r$ and $v$ slightly smaller so that
$r<\e/4$. Fix $\d<r/100$ and find a path in $U_k^{\e/2}$ of length
between the numbers $\mathrm{diam}^{int} (U_k^{\e/2})$ and
$\mathrm{diam}^{int} (U_k^{\e/2})+\d$. This path is almost length
minimizing with the error $\le\d$. Hence, one can find
$N=[\mathrm{diam}^{int} (U_k^{\e/2})/3r]$ points on the path such
that $r$-balls centered at the points are disjoint. Thus, $V\ge
\mathrm{vol}(U_k^{\e})\ge Nv$ and we get a uniform bound on
$\mathrm{diam}^{int} (U_k^{\e/2})$.

Conversely, let us prove (ii)$\Rightarrow$(i). Fix $\d <\e/10$.
First, show that there is a uniform lower bound for
$\mathrm{vol}(B(x,\d))$ for any $x\in U_k^{\e/2}$. Take an
arbitrary point $x\in U_k^\e/2$. Since $\mathrm{diam}^{int}
(U_k)\le D$, there is a sequence of points $x_i\in U_k$, $i=0,
\dots, N$ where $N=[D/\d]+1$ that starts at $x_0$, ends at $x_N=x$
and satisfies $d(x_i,x_{i+1})\le\d$. Let $v_n(r,H)$ denote the
volume of the ball of radius $r$ in a complete simply connected
$n$-dimensional space of constant sectional curvature $=H$. Using
induction on $i$ and the relative volume comparison, one can show
that for every $i$ $$\mathrm{vol}(B_\d(p_i))\ge
\left(\frac{v_n(\d,H)}{v_n(2\d,H)}\right)^iv.$$ In particular,
there is a uniform lower bound for $\mathrm{vol}(B(x,\d))$.

Now fix a finite covering of $U_k^{\e/2}$ by $\d$-balls.
As before the relative volume comparison
gives a uniform upper bound $N_{loc}(\d)$
on the multiplicities of intersections in this covering
(the argument involves only small balls so it works because
the balls are far enough from the boundary.)

By the absolute volume comparison, the volume of each $\d$-ball
is uniformly bounded above, hence a bound $\mathrm{vol}(U_k^{\e})\le V$
would follow from a bound on the number of balls in the covering.
Set $r_j=j\d$, $j=0,\dots, m$ with $m=[D/\d]+1$.
Let $N_j$ be the number of balls in the covering whose
centers are in the $r_j$-ball around $x_0$
(as before the ball is taken with respect to the induced
Riemannian metric on $U_k^{\e/2}$).
Since multiplicities are bounded by $N_\mathrm{loc}$,
for each $j$ we have that
$N_{j+1}\le N_j+N_j^{N_{loc}(\d)}$. This
gives a uniform bound on the number of balls in the covering, and hence on
$\mathrm{vol}(U_k^{\e})$.

\end{rmk}

\section{Invariants of maps}
\label{Invariants of maps}

\begin{defn}
Let $B$ be a topological space and $S(B)$ be a set.  Given a
smooth manifold $N$, denote by $C(B,N)$ the set of all continuous
maps from $B$ to $N$. Suppose that for any smooth $N$ we have a
map $\iota \co C(B,N)\to S(B)$. Then we call $\iota$ an
$S(B)${\it -valued invariant of maps of $B$} if the two following conditions
hold:

$\mathrm{(1)}$ Homotopic maps $f_1:B\to N$ and
$f_2\co B\to N$ have the same invariant.

$\mathrm{(2)}$ Let $h\co  N\to L$ be a
homeomorphism of $N$ onto an open subset of $L$.
Then, for any continuous map $f\co B\to N$,
the maps $f\co B\to N$ and $h\circ f\co B\to L$
have the same invariant.
\end{defn}

There is a variation of this definition for maps into
oriented manifolds.
Namely, we require that the target manifold is oriented and
the homeomorphism $h$ preserves orientation.
In that case we say that $\iota$ is an
{\it invariant of maps into oriented manifolds}.

\begin{ex}\bf \!(Pontrjagin classes)\rm\
\label{Rational Pontrjagin classes} As usual the total (rational)
Pontrjagin class of a bundle $\xi$ is denoted by $p(\xi)$. Given a
continuous map of smooth manifolds $f\co B\to N$, set $p(f)$ to be
the solution of $f^\ast p(TN)=p(TB)\cup p(f)$. (A total Pontrjagin
class is a unit so there exists a unique solution.) The fact that
$p(f)$ is an $H^\ast(B,\mathbb Q)$-valued invariant follows from
topological invariance of rational Pontrjagin classes~\cite{Nov}.
In case $f$ is a smooth immersion, $p(f)$ is the the total
Pontrjagin class of the normal bundle to $f$. Finally, note that
Stiefel-Whitney classes are preserved by homeomorphism~\cite{Spn}
hence they also give rise to invariants of maps.
\end{ex}

\begin{rmk}
The isomorphism class of the pullback of the tangent bundle to $N$
under $f$ would be an invariant (for paracompact $B$) if we only
require that invariants are preserved by diffeomorphisms. In
general, homeomorphisms do not preserve tangent bundles. However,
tangent bundle (and, in fact, any vector bundle over a finite cell
complex) is recovered up to finitely many possibilities by the
total Pontrjagin class and the Euler class of its orientable ($1$
or $2$-fold) cover (see~\cite{Bel} for a proof of this folklore
result).
\end{rmk}

\begin{ex}\bf \!(Intersection number in oriented $n$-manifolds.)\rm\
\label{intersection number}
Assume $B$ is a compact space
and fix two homology classes
$\alpha\in H_{m}(B,\mathbb Q)$ and $\beta\in H_{n-m}(B,\mathbb Q)$.
(Unless stated otherwise we always use singular (co)homology
with rational coefficients.)
Let $f\co B\to N$ be a continuous map of a
compact topological space $B$ into
an oriented $n$-manifold $N$.
Set $I_{n,\alpha,\beta}(f)$ to be the
intersection number
of $f_*\alpha$ and $f_*\beta$ in $N$.
It is easy to see that $I_{n,\alpha,\beta}$
is an $\mathbb Q$-valued invariant of maps into
oriented manifolds.
\end{ex}

\begin{ex}\bf \!(Generalized Euler class)\rm\
\label{euler class}
Let $B$ be a closed oriented $m$-manifold
and let $f\co B\to N$ be a map of $B$ into an oriented
$n$-manifold $N$.
Define the rational Euler class $\chi (f)$ by
requiring that $\langle \chi (f), \alpha\rangle=I_{n,\alpha,[M]}$.
This is a $H^{n-m}(B,\mathbb Q)$-valued invariant
for maps into oriented manifolds.
If $f$ is a smooth embedding, $\chi (f)$ is the Euler class of
the normal bundle $\nu_f$.
Note that when the orientation is changed on $B$ or $N$,
the invariant $\chi (f)$ may change sign.

More generally, if $B$ and $N$ are not assumed to be orientable
one can define a (generalized) Euler class of a continuous map as follows.

Recall that a smooth manifold $L$ is orientable  iff the first
Stiefel-Whitney class $w_1(TL)\in H^1(L,\mathbb Z/2\mathbb Z)$
vanishes. Note that $$H^1(L,\mathbb Z/2\mathbb Z)\cong
\mathrm{Hom}(H_1(L),\mathbb Z/2\mathbb Z)\cong
\mathrm{Hom}(\pi_1(L),\mathbb Z/2\mathbb Z),$$ so elements of
$H^1(L,\mathbb Z/2\mathbb Z)$ correspond to subgroups of index
$\le 2$ in $\pi_1(L)$ which are the kernels of homomorphisms in
$\mathrm{Hom}(\pi_1(L),\mathbb Z/2\mathbb Z)$.

Let $K_f$ be the intersection of the subgroups corresponding to
$w_1(B)$ and $f^\ast w_1(N)$. Let $\widetilde B\to B$ be a
covering associated to $K_f$ and let $\widetilde N\to N$ be a
covering associated with $f_\ast (K_f)$. Then $f$ lifts to a map
$\tilde f\co\widetilde B\to\widetilde N$ of orientable manifolds.
Define the {\it generalized Euler class} $\tilde\chi(f)$ as a pair
$(K_f, \pm \chi (\tilde f))$ (note that $\chi (\tilde f)$ depends
on the choice of orientations in $\widetilde B$ and $\widetilde
N$, so it is only well-defined up to sign). It is easy to see that
$\tilde\chi(f)$ is an invariant because homotopies and
homeomorphisms lift to covering spaces, and because
Stiefel-Whitney classes are topological invariants~\cite{Spn}.

Thus, for our purposes,  $\tilde\chi(f)$ is a regular covering
$\widetilde B\to B$ and two cohomology classes
$\chi (\tilde f)$, $-\chi (\tilde f)$ in
$H^{n-m}(\widetilde B, \mathbb Q)$.
For a map of orientable manifolds $f\co B\to M$,
$\tilde\chi(f)=(\pi_1(B), \pm\chi(f))$ so, up to sign,
$\tilde\chi(f)$ generalizes $\chi(f)$.

If $f$ is a smooth embedding of nonorientable manifolds
with orientable normal bundle $\nu_f$, then
the Euler class of $\nu_f$ is taken to $\pm\chi(\tilde f)$
by the map $H^{n-m}(B,\mathbb Q)\to H^{n-m}(\widetilde B,\mathbb Q)$
induced by the covering.
\end{ex}

\begin{prop} \label{pe determine bundles}
Let $e_k\co B\to N_k$
be a sequence of smooth embedding of a closed manifold $B$
into manifolds $N_k$. Assume that invariants $p(e_k)$ and
$\tilde\chi(e_k)$ are independent of $k$.
Then the set of isomorphism classes of normal bundles
$\nu_{e_k }$ is finite.
\end{prop}
\begin{proof} It is well-known to experts that
a vector bundle over a finite cell complex is recovered up to
finitely many possibilities by the total Pontrjagin class and the
Euler class of its orientable ($1$ or $2$-fold) cover
(see~\cite{Bel} for a proof). We are now going to reduce to this
result.

In what follows we use the  notations of~\ref{euler class}. Since
$\tilde\chi(e_k)=(K_{e_k}, \pm\chi (\tilde e_k))$ is independent
of $k$, there is a covering $\widetilde B\to B$ associated with
$K_{e_k}\le\pi_1(B)$ and, for each $k$, a covering $\widetilde
N_k\to N_k$ associated with $e_{k\ast}(K_{e_k})$. The embedding
$e_k$ lifts to an embedding $\tilde e_k\co\widetilde
B\to\widetilde N_k$ of orientable manifolds.

Using that $H^1(B,\mathbb Z/2\mathbb Z)$ is a finite group, we can
partition $\nu_{e_k }$ into finitely many subsequences each
having the same first Stiefel-Whitney class. It suffices
to show that any such subsequence falls into finitely many
isomorphism classes, so we can assume that $w_1(\nu_{e_k })$, and hence
$e_k^\ast w_1(TN_k)=w_1(\nu_{e_k })+w_1(TB)$, is independent of $k$.
Let $\overline B\to B$ be a covering associated
to the subgroup of $\pi_1(B)$ that
corresponds to $w_1(\nu_{e_k })$.
This subgroup lies in $K_{e_k}$ so $\overline B$ is an intermediate
covering space between $\widetilde B$ and $B$, that is,
we have coverings $\tilde q\co\widetilde B\to\overline B$
and $\bar q\co\overline B\to B$.
Also let $\overline N_k\to N_k$ be a covering associated to
$e_{k\ast}\bar q_\ast(\pi_1(\overline B))$.

The embedding $e_k $ lifts to an embedding
$\bar e_k \co\overline B\to\overline N$. Now the normal
bundle $\nu_{\bar e_k}$ is orientable so its Euler class is well-defined
(up to sign since there is no canonical choice of orientations).
Note that $\tilde q^\ast$ takes the Euler class of $\nu_{\bar e_k}$ to
$\pm\chi (\tilde e_k)$.
It is a general fact that finite covers induce
injective maps in rational cohomology.
(The point is that the transfer map goes the other way,
and precomposing the transfer
with the homomorphism induced by the covering is
multiplication by the order of the covering.)
Also, the Pontrjagin class of $\nu_{\bar e_k}$ is
$\bar q^\ast$-image of $p(e_k )$.

Thus, given $\chi (\tilde e_k)$ and $p(e_k )$,
one can uniquely recover the (rational) Euler and Pontrjagin class of
$\nu_{\bar e_k}$.
As we mentioned above these classes determine
$\nu_{e_k }$ up to finitely many possibilities.
Therefore, $\nu_{e_k }$ are determined up to finitely many possibilities by
$\tilde\chi (e_k)$ and $p(e_k )$ as desired.
\end{proof}

\begin{rmk} \label{pe determine bundles rmk}
The above proof actually gives a slightly more general result
which will be useful in our applications. Namely, instead of
assuming that $e_k $'s are smooth embeddings, it suffices to
assume that each $e_k $ is a topological embedding such that $e_k
(B)$ is a smooth submanifold of $N_k$. Then $e_k (B)$ has a normal
bundle in $N_k$ whose pullback via $e_k $ is still denoted
$\nu_{e_k }$.
\end{rmk}

\section{Main technical results}
\label{Main observation}

\begin{prop}\label{main observation}
Let $f_k\co M\to N_k$ be a sequence of continuous maps of a closed
Riemannian manifold $M$ into (possibly incomplete) Riemannian
$n$-manifolds. Assume that, for each $k$, there exists a compact
domain $U_k\supset f_k(M)$ such that $\{U_k\}$ has uniformly
bounded geometry. Assume that either \vspace{2pt}\newline
\hspace*{10pt}\rm (i)~\it$ \{f_k\}$ is almost equicontinuous, or
\vspace{2pt}\newline \hspace*{10pt}\rm (ii)~\it $f_k$ is a
homotopy equivalence with a homotopy inverse $g_k\co N_k\to M$
such\newline \hspace*{10pt}that $\{g_k\}$ is almost
equicontinuous. \vspace{2pt}\newline Then, for any $S(M)$-valued
invariant of maps $\iota$, the subset $\{\iota(f_k)\}$ of $S(M)$
is finite.
\end{prop}
\begin{proof} Since $\{U_k\}$ has uniformly bounded geometry,
there exists a metric space $V$, and homeomorphisms
$h_k\co V\to V_k$ of $V$ onto compact domains $V_k\supset U_k$ such that
both $\{h_k\}$ and $\{h_k^{-1}\}$ are
almost equicontinuous.

If (i) holds, then $\{h_k^{-1}\circ f_k\}$ is an almost
equicontinuous sequence of maps from $M$ into $V$.
Thus,~\ref{AAthm} implies that the maps $h_k^{-1}\circ f_k$'s fall
into finitely many homotopy classes. Now we are done by definition
of an invariant since $h_k^{-1}$'s are homeomorphisms.

If (ii) holds, then $\{g_k\circ h_k\}$ is an almost equicontinuous
sequence of maps from $V$ into $M$.
Again, by~\ref{AAthm}
there are only finitely many homotopy classes of maps among
$g_k\circ h_k$.
It suffices to show that, whenever $g_k\circ h_k$ is homotopic to $g_m\circ
h_m$, the maps $f_k$ and $f_m$ have the same invariants.
Let $G\co V\times[0,1]\to M$ be a homotopy that connects
$g_k\circ h_k$ and $g_m\circ h_m$.
The homotopy $F\co M\times [0,1]\to V$ defined as
$F(x,t)=G(h_k^{-1}(f_k(x)),t)$ connects
$g_m\circ h_m\circ h_k^{-1}\circ f_k$ with
$g_k\circ h_k\circ h_k^{-1}\circ f_k=g_k\circ f_k\sim \mathrm{id}_M$.
Thus, $f_m$ is homotopic to
$$f_m\circ g_m\circ h_m\circ h_k^{-1}\circ f_k\sim
\mathrm{id}_{N_m}\circ h_m\circ h_k^{-1}\circ f_k\sim
h_m\circ h_k^{-1}\circ f_k.$$
Since $h_m\circ h_k^{-1}$ is a homeomorphism,
$f_k$ and $f_m$ have the same invariants as desired.
\end{proof}

\begin{rmk}
In the above theorem $M$ can be chosen as in~\ref{rmk after AA
thm}. Note that the spaces mentioned in~\ref{rmk after AA thm} are
homotopy equivalent to finite cell complexes~\cite{West}, in
particular, characteristic classes determine a vector bundle over
such a space up to finitely many possibilities.

Also, instead of assuming $\{g_k\}$ is almost equicontinuous, it
is enough to assume that $g_k|_{V_k}$ is almost equicontinuous.

There is a version of the theorem for invariants of maps into
oriented manifolds. First of all, by pulling back the orientation
from $V$, one can {\it define} orientations on $N_k$ so that
$h_k$ preserve orientations. In general, change of orientation on $N$
may lead to an unknown change of an invariant of a map into $N$.
However, if $\iota=I_{n,\alpha,\beta}$, then change of orientation on $N$
may only lead to the sign change for the intersection number.
Thus, for $\iota=I_{n,\alpha,\beta}$, the above theorem holds.
\end{rmk}

\begin{cor}
\label{main observ: fixed base}
Let $M$ be a closed Riemannian manifold
and let $e_k\co M\to N_k$ is a sequence of topological embeddings
of $M$ into Riemannian $n$-manifolds $N_k$ such that $e_k(M)\subset N_k$
is a smooth submanifold. Assume that, for each $k$,
$e_k$ is homotopic to $f_k\co M\to N_k$ and
there exists a compact domain $U_k\supset f_k(M)$
such that $\{U_k\}$ has uniformly bounded geometry.
Assume that either
\vspace{2pt}\newline
\hspace*{10pt}\rm (i)~\it$ \{f_k\}$ is almost equicontinuous, or
\vspace{2pt}\newline
\hspace*{10pt}\rm (ii)~\it $f_k$ is a homotopy equivalence with a homotopy
inverse
$g_k\co N_k\to M$ such\newline
\hspace*{10pt}that $\{g_k\}$
is almost equicontinuous.
\vspace{2pt}\newline
Then the set of isomorphism classes of normal bundles
$\nu_{e_k}$ is finite.
\end{cor}
\begin{proof}
For any invariant $\iota$, $\iota(f_k)=\iota(e_k)$.
In particular this is true for the rational Pontrjagin class
and generalized Euler class.
The result now follows from~\ref{main observation}
combined with~\ref{pe determine bundles},~\ref{pe determine bundles rmk}.
\end{proof}

\begin{rmk} Note that~\ref{main observ: fixed base} also holds when
$e_k$'s are only immersions provided $\dim(N_k)-\dim(M)$ is either
odd or $>\dim(M)$. Indeed, under these codimension assumptions the
rational Euler class of $\nu_{e_k}$ vanishes while the total
Pontrjagin class of $\nu_{e_k}$ is equal to $p(e_k)$.
\end{rmk}

\begin{cor}
\label{main observ: variable base}
Let $e_k\co M_k\to N_k$ be a sequence of topological embeddings
of closed Riemannian manifolds $M_k$ into Riemannian $n$-manifolds $N_k$
such that $e_k(M_k)\subset N_k$ is a smooth submanifold
and $\{M_k\}$ has uniformly bounded geometry.
Assume that
$e_k$ is homotopic to $f_k\co M_k\to N_k$ and
there exists a compact domain $U_k\supset f_k(M_k)$
such that $\{U_k\}$ has uniformly bounded geometry.
Assume that either
\vspace{2pt}\newline
\hspace*{10pt}\rm (i)~\it$ \{f_k\}$ is almost equicontinuous, or
\vspace{2pt}\newline
\hspace*{10pt}\rm (ii)~\it $f_k$ is a homotopy equivalence with a homotopy
inverse
$g_k\co N_k\to M_k$ \newline
\hspace*{10pt}such that $\{g_k\}$
is almost equicontinuous.
\vspace{2pt}\newline
Then the set of topological equivalence classes of normal bundles
$\nu_{e_k}$ is finite.
\end{cor}
\begin{proof}
Since $M_k$ has bounded geometry, there exists $M$
and homeomorphisms $h_k\co M\to M_k$ such that
both $\{h_k\}$ and  $\{h_k^{-1}\}$ are almost equicontinuous.
Note that $e_k(h_k(M))=e_k(M_k)$ is a smooth submanifold of $N_k$.
If $\{f_k\}$ is  almost equicontinuous, then so is $\{f_k\circ h_k\}$.
Similarly, if $\{g_k\}$ is almost equicontinuous, then so is
$\{h_k^{-1}\circ g_k\}$.
Thus,~\ref{main observ: fixed base} implies that
the set of isomorphism classes of normal bundles
$\nu_{e_k\circ h_k}$ is finite.
In particular, the set of topological equivalence classes of normal bundles
$\nu_{e_k}$ is finite.
\end{proof}

\begin{rmk} For future applications we note that
if $h_k$'s are diffeomorphisms, then the conclusion of~\ref{main
observ: variable base} can clearly be improved to ``the set of
smooth equivalence classes of normal bundles $\nu_{e_k}$ is
finite.''
\end{rmk}

\section{Geometric applications}
\label{sec: appl}

This section contains proofs of the various finiteness theorems
that follow from section~\ref{Main observation}.

\begin{cor}\label{appl: cor of fixed base}
Let $e_{\alpha}\co M\to N_{\alpha}$ be an almost equicontinuous
family of smooth embeddings of a closed Riemannian manifold $M$
into complete Riemannian $n$-manifolds $N_{\alpha}$ with
$\sec(N_{\alpha})\ge -1$. Assume that for each $\alpha$ there is a
point $p_{\alpha}\in N_{\alpha}$ such that
$\mathrm{vol}(B(p_{\alpha},1))$ is uniformly bounded below and
$dist_{N_{\alpha}}(p_{\alpha}, e_{\alpha}(M))$ is uniformly
bounded above. Then the set of isomorphism classes of normal
bundles $\nu_{e_{\alpha}}$ is finite.
\end{cor}
\begin{proof}
Since $\{e_\alpha\}$ is almost equicontinuous,
$\mathrm{diam}(e_\alpha(M))$ is uniformly bounded above. The
result now follows from~\ref{main observ: fixed base} and
section~\ref{sec: local}.
\end{proof}

\begin{proof}[Proof of~\ref{intro: fixed soul}]
Since $\mathrm{diam}(f(S))\le D$ we can find a compact domain
$U\supset f(S)$ with $\mathrm{diam}(U)\le 2D$.
By results of the section~\ref{sec: local},
any family of such domains $U$
has bounded geometry, hence the conclusion follows
from~\ref{main observ: fixed base}.
\end{proof}

\begin{proof}[Proof of~\ref{intro: variable soul}]
Let $N_{\alpha}$ be a family of nonnegatively curved manifolds
satisfying conditions of~\ref{intro: variable soul}. For any
$\alpha$ let $S_{\alpha}\subset N_{\alpha}$ be a soul of
$N_{\alpha}$. First,we show that $\{S_{\alpha}\}$ has uniformly
bounded geometry. By assumption $S_{\alpha}$ has lower volume
bound. Lower sectional curvature bound follows because souls are
totally geodesic. Since $\mathrm{diam}(f_{\alpha}(S_{\alpha}))\le
D$ and there is a distance-nonincreasing retraction of
$r_{\alpha}\co N_{\alpha}\to S_{\alpha}$~\cite{Sha}, the diameter
of $r_{\alpha}(f_{\alpha}(S_{\alpha}))$ is at most $D$. The map
$r_{\alpha}\circ f_{\alpha}\co S_{\alpha}\to S_{\alpha}$ is a
homotopy equivalence, in particular, it has nonzero degree, hence
it is onto. We conclude that $\mathrm{diam}(S_{\alpha})\le D$.
Thus, $\{S_{\alpha}\}$ has uniformly bounded geometry.

Since $\mathrm{diam}(f_{\alpha}(S_{\alpha}))\le D$ we can find
compact domains $U_{\alpha}\supset f_{\alpha}(S_{\alpha})$ with
$\mathrm{diam}(U_{\alpha})\le 2D$. Again, $\{U_{\alpha}\}$ has
uniformly bounded geometry, and the conclusion follows
from~\ref{main observ: variable base}.
\end{proof}

We now prove a theorem that, in particular,
implies~\ref{intro: no collapse on soul}.

\begin{thm}\label{appl: lower sec, tot geod}
Given $n$, $K$, $D$, $r$, $v$, there is a finite
collection of vector bundles such that for any
totally geodesic embedding of a closed Riemannian manifold
$M$ into a complete Riemannian manifold $N$, the normal bundle
of $M$ is topologically equivalent to a bundle
of the collection provided
$\mathrm{diam}(M)\le D$, $\sec(N)\ge -1$, and
there exist positive $r$, $v$, and a point $p\in e(M)$ such that
$\mathrm{vol}_{N}B(p,r)\ge v$.
\end{thm}

\begin{proof} Start with an arbitrary family of
totally geodesic embeddings $e_{\alpha}\co
M_{\alpha}\hookrightarrow N_{\alpha}$ as above. First, we show
that $\{M_{\alpha}\}$ has uniformly bounded geometry where
$M_{\alpha}$ is equipped with the induced Riemannian metric. By a
result of Karcher-Heinze~\cite{HK}
$\mathrm{vol}_{N_{\alpha}}B(p_{\alpha},r)\ge v$ implies a lower
volume bound on $M_{\alpha}$. Since $M_{\alpha}$ is totally
geodesic $\sec(M_{\alpha})\ge -1$, and by assumption
$\mathrm{diam}(M_{\alpha})\le D$. Thus, Perelman's stability
theorem implies that $\{M_{\alpha}\}$ has uniformly bounded
geometry (see~\ref{local: perelman})

Note that $\mathrm{diam}(e_{\alpha}(M_{\alpha}))\le
\mathrm{diam}(M_{\alpha})\le D$, hence,~\ref{local: perelman} implies
that there is a compact domain $W_{\alpha}\supset
e_{\alpha}(M_{\alpha})$ such that $\{W_{\alpha}\}$ has uniformly
bounded geometry. The result now follows from~\ref{main observ:
variable base}(ii) because totally geodesic embeddings
$\{e_{\alpha}\}$ are $1$-Lipschitz, in particular,
$\{e_{\alpha}\}$ is equicontinuous.
\end{proof}

The same proof gives the following.

\begin{thm}\label{appl: lower ric, tot geod}
Given positive $n$, $H$, $D$, $i_0$, and $\e$, there is a finite
collection of vector bundles such that for any
totally geodesic embeddings of a closed Riemannian manifold
$M$ into a complete Riemannian manifold $N$, the normal bundle
of $M$ is topologically equivalent to a bundle
of the collection provided
$\mathrm{diam}(M)\le D$, $\mathrm{Ric}(N)\ge -1$, and $\mathrm{inj}(x)\ge
i_0$ for any
$x$ in the $\e$-neighborhood of image of $M$.{\hfill\qed}
\end{thm}

\begin{rmk} There are obvious ``fixed base'' modifications
of~\ref{appl: lower sec, tot geod}
and~\ref{appl: lower ric, tot geod}.
\end{rmk}

\begin{cor} \label{appl: euler class and no collapse on soul}
Given positive $D$, $r$, $v$, $n$, and a closed manifold $M$ with
$\oplus_{i>0}H^{4i}(M,\mathbb Q)=0$,
there exists a finite collection of vector bundles over $M$ such that,
for any open complete Riemannian $n$-manifold $N$ with $\sec(N)\ge 0$
and a soul $S\subset N$, the normal bundle to $S$ is
topologically equivalent to a bundle of the collection provided
$S$ is homeomorphic to $M$, and
the inclusion $S\hookrightarrow N$
is homotopic to a map $f$ such that
$\mathrm{diam}(f(S))\le D$ and $\mathrm{vol}_{N}(B(p, r))\ge v$
for some $p\in f(S)$.
\end{cor}
\begin{proof}
First note that, up to topological equivalence, only
finitely many of the bundles $\nu_S$ can have zero
Euler class. (Otherwise, there is a sequence of pairwise topologically
inequivalent bundles $\nu_{S_k}$ with zero Euler class.
Use homeomorphisms $M\to S_k$ to pull the bundles back
to $M$. These pullback bundles clearly have zero
Euler class as well as zero rational Pontrjagin classes since
$\oplus_{i>0}H^{4i}(M,\mathbb Q)=0$. Thus the bundles
belong to finitely many isomorphism classes which implies that
$\nu_{S_k}$ belong to finitely many topological equivalence
classes.)
Now if the Euler class
of the normal bundle to $S$ is nonzero, then
$f(S)\cap S\neq\emptyset$, hence the distance from
$p$ to $S$ is $\le D$,
and the result follows from~\ref{appl: lower sec, tot geod}.
\end{proof}

\begin{proof}[Proof of~\ref{intro: fixed soul, normal curv}]
Let $N^n$ and $p\in S\subset N$ be
chosen to satisfy the assumptions. According to~\ref{intro: fixed soul}
we only have to show
is that under our  assumptions we have a uniform lower  bound on
$\mathrm{vol} B(p,r)$. Let $r_0=\min \{r/2, \pi/(2\sqrt{K})\}$.  Let
$l=codim
S-1$ and  $(S^l_K,g_{can})$ be a round sphere of constant
curvature $K$ and $\bar{p}$ be  any point on this sphere. Consider
the exponential map $\exp_K\co T_{\bar{p}}S^l_K\rightarrow S^l_K$.
Denote by $v(l,K,t)$ the volume of the ball of radius $t$ centered
at $\bar{p}$.

First of all, notice that by the triangle inequality
$B(p,r)$ contains a tubular neighborhood $U(p,r_0)$  consisting of
all points $x\in N$ such that $d(x,S)\le r_0$ and $d(p,Sh(x))\le
r_0$. Here $Sh$ stands for the Sharafutdinov retraction $Sh\co
N\rightarrow S$. For any $x\in S$ denote 
$$B^\perp(x,t)=\{y\in N| d(y,S)\le t \text{ and } Sh(y)=x\}.$$ 
Since Sharafutdinov
retraction is a $C^1$-Riemannian submersion~\cite{Per2} we can
apply Fubini's Theorem to see that
\begin{equation}\label{fubini}
vol U(p,r_0)=\int_{B_S(p,r_0)}\mathrm{vol} B^\perp(x,r_0)\mathrm{dvol}(x)
\end{equation}
Here $B_S(p,r_0)$ stands for the ball of radius $r_0$ around $p$
in $S$. It suffices to show that for each $x\in B_S(p,r_0)$ we
have $\mathrm{vol} B^\perp(x,r_0)\ge v(l,K,r_0)$.
(Indeed, it would imply that
$\mathrm{vol}(B(p,r))\ge \mathrm{vol}(U(p,r_0))\ge
\mathrm{vol}(B_S(p,r_0))\cdot v(l,K,r_0)$.
Finally, by volume comparison, $\mathrm{vol}(B_S(p,r_0))$ is bounded below
in terms of
$D$ and $v$ and we are done.)

Fix an $x\in B_S(p,r_0)$ and consider the normal exponential map
$\exp_x^\perp\co T_x^\perp S\rightarrow N$.  It follows from~\cite{Per2}
that this map sends the ball $B_{T_x^\perp}(0,r_0)$
onto $B^\perp(x,r_0)$. Choose a linear isometry between
$T_x^\perp$ and $T_{\bar{p}}S^l_K$ and use it to equip
$B_{T_x^\perp}(0,r_0)$ with the metric $ g_K$ of constant
curvature $K$. Let $g_x$ be the induced Riemannian metric on the
Sharafutdinov fiber over $x$. To finish the proof it is enough to
establish the following lemma saying that "reverse Toponogov
comparison" holds on $B^\perp(x,r_0)$.

\begin{lem}
The surjection $\exp_x^\perp\co (B_{T_x^\perp }(0,r_0),g_K)\rightarrow
(B^\perp(x,r_0),g_x)$ is a distance nondecreasing diffeomorphism.
\end{lem}
Let $v$ be a unit vector in $T_x^\perp$ and $\gamma(t)= \exp(tv)$
be the normal geodesic in direction $v$.
%
%
We now show that, for any $t\le r_0$ and any $X\in T_{\gamma(t)}$ with
$|X|=1$ and $\langle X,\gamma'(t)\rangle$, we have that $K(X,\gamma'(t))\le
K$.
Write $X=X^h+X^v$ as a sum of its horizontal and vertical
components. Then

$K(X,\gamma'(t))=
\langle R(\gamma'(t),X^h+X^v)\gamma'(t),X^h+X^v\rangle =
\langle R(\gamma'(t),X^v)\gamma'(t),X^v\rangle +
\langle R(\gamma'(t),X^h)\gamma'(t),X^h\rangle +
\langle R(\gamma'(t),X^h)\gamma'(t),X^v\rangle +
\langle R(\gamma'(t),X^v)\gamma'(t),X^h\rangle .$

The first term in the right hand side is $\le K$ by assumption and
also because $|X^v|\le |X|=1$ .  By~\cite{Per2}
$R(\gamma'(t),X^h)\gamma'(t)=0$ and therefore
$$\langle R(\gamma'(t),X^h)\gamma'(t),X^v\rangle= \langle
R(\gamma'(t),X^h)\gamma'(t),X^h\rangle=0.$$
By the symmetry of the
curvature tensor the forth term is equal to the third one and
hence is also equal to $0$. Thus $K(X,\gamma'(t)= \langle
R(\gamma'(t),X^h)\gamma'(t),X^h\rangle\le K$.

Now since $r_0<\pi/\sqrt{K}$ and all the two planes along $\gamma(t)$
containing $\gamma'(t)$ have curvature $\le K$ we can
apply the Rauch comparison theorem  to conclude that the
differential of $\exp_x^\perp$ does not decrease the lengths of
tangent vectors and thus
$\exp_x^\perp\co (B_{T_x^\perp}(0,r_0),g_K)\rightarrow (B^\perp(x,r_0),g_x)$
is a local diffeomorphism that does not decrease lengths of curves.
It remains to show that this map is $1-1$.

Suppose not. Then the injectivity radius $r_x$ of $B^\perp(x,r_0)$
at $x$ is strictly less than $r_0$. Let $v,u\in T_x^\perp$ be such
that $|u|= |v|= r_x$ and $\exp(v)= \exp(u)$. Denote $q= \exp(v)=
\exp(u)$. Notice that geodesics $\gamma_v(t)= \exp(tv)$ and
$\gamma_u(t)= \exp(tu)\co [0,1]\rightarrow N$ connecting $x$ and
$q$ are obviously distance minimizing. By~\cite[Lemma~5.6.5]{CE}
 these geodesics must form a geodesic loop at $x$ (i.e
$\gamma_u'(1)= -\gamma_v'(1)$). This is impossible since according
to~\cite{CG} the soul $S$ is totally convex and $x$ lies in $S$.
\end{proof}

\begin{cor}\label{appl: fixed soul, normal curv}
Given positive $D$, $n$, $v$, $r$, and $K$, there exists a finite family of
vector bundles such that, for any complete open Riemannian $n$-manifold
$N$ with $\sec(N)\ge 0$ and a soul $S\subset N$,
the normal bundle to $S$ is topologically
equivalent to a bundle of the collection provided
$\mathrm{diam}(S)\le D$, $\mathrm{vol}(S)\ge v$ and there is a point
$p\in S$ such that all the vertical curvatures at the points of
$B(p,r)$ are bounded above by $K$.
\end{cor}
\begin{proof} The proof of~\ref{intro: fixed soul, normal curv}
gives a uniform lower bound on $\mathrm{vol} B(p,r)$ so the result
follows from~\ref{intro: variable soul}.
\end{proof}

\begin{cor} \label{appl: fixed nonpositive base}
Given positive $D$, $r$, $v$, $n$, $K$, and a closed
Riemannian manifold $M$ with $\sec(M)\in [-1,0]$,
there exists a finite collection of vector bundles over $M$ such that,
for any totally geodesic embedding $e\co M\to N$ of $M$
into an open complete Riemannian $n$-manifold $N$ with $\sec(N)\le 0$,
the normal bundle $\nu_e$ is isomorphic to
a bundle of the collection provided
$\mathrm{vol}(M)\ge v$ and $e$
is  a homotopy equivalence homotopic to a map $f$ with
$\mathrm{diam}(f(M))\le D$ such that the sectional curvature at any point of
the $r$-neighborhood of $f(M)$ is $\ge -K$.
\end{cor}
\begin{proof}
Start with an arbitrary family of
totally geodesic embeddings
$e_{\alpha}\co M\hookrightarrow N_{\alpha}$ as above.
First, note that for any $\alpha$
 the injectivity radius of $N_{\alpha}$ satisfies
$\mathrm{inj}(N_{\alpha})\ge
\mathrm{inj}(M)=\mathrm{inj}(e_{\alpha}(M))$. If not, there is a
point $p\in N_{\alpha}$ with
$\mathrm{inj}_{N_{\alpha}}(p)<\mathrm{inj}(M)$. Since $N_{\alpha}$
and $M$ has nonpositive sectional curvatures, the injectivity
radius at any point is half the length of the shortest geodesic
loop at this point~\cite[Lemma~5.6.5]{CE}. Take a geodesic loop at
$p$ of length $<\mathrm{inj}(M)/2$ and project it to
$e_{\alpha}(M)$ by the closest point retraction $r_{\alpha}\co
N_{\alpha}\to e_{\alpha}(M)$~\cite{BGS}. The retraction is a
distance-nonincreasing homotopy equivalence. (In fact, since
$\sec(N_{\alpha})\le 0$, the normal exponential map identifies
$N_{\alpha}$ with the normal bundle to $e_{\alpha}(M)$ where
$r_{\alpha}$ corresponds to the bundle projection.) Thus, we get a
homotopically nontrivial curve of length $<\mathrm{inj}(M)/2$ in
$e_{\alpha}(M)$. Since $e_{\alpha}$ is an isometric embedding it
preserves lengths of curves. Therefore, we obtain a homotopically
nontrivial curve of length $<\mathrm{inj}(M)/2$ in $M$ which is
impossible because loops of length $<\mathrm{inj}(M)/2$ lift to
loops in the universal cover.

Find compact domains $U_{\alpha}\supset f_{\alpha}(M)$ that lie in
the $r$-neighborhood of $f_{\alpha}(M)$. Since
$\mathrm{diam}(f_{\alpha}(M))\le D$, the diameter of $U_{\alpha}$
is $\le D+2r$. Also $\sec(U_{\alpha})\in [-K,0]$ and
$\mathrm{inj}(U_{\alpha})\ge \mathrm{inj}(M)$,
hence~\ref{local: anderson-cheeger}
implies that $\{U_{\alpha}\}$ has bounded geometry, and the
conclusion follows from~\ref{main observ: fixed base}.
\end{proof}

\begin{proof}[Proof of~\ref{intro: variable nonpositive base}]
Since $\sec(M_{\alpha})$ is bounded in absolute value, lower bound on
volume implies a lower bound on the injectivity radius. Thus, the
result follows from~\ref{main observ: variable base} exactly as in
the proof of~\ref{appl: fixed nonpositive base}.
\end{proof}
\begin{rmk} In the statement of~\ref{intro: variable nonpositive base}
one can replace ``$\mathrm{vol}(M)\ge v$'' by ``$\pi_1(M)$ has no
normal virtually abelian subgroups'', or equivalently, by the
universal cover of $M$ has no Euclidean de Rham factor (see
e.g.~\cite[pp395--396]{Fu}).
\end{rmk}

\begin{rmk} Given a closed Riemannian manifold $M$, there are only
finitely many isomorphism classes of normal bundles of isometric
immersions $f\co M\to N$ into Riemannian manifolds $N$ such that
$|\sec(N)|$ and the second fundamental form of $f$ are uniformly
bounded. Indeed, these bounds imply a uniform bound on the
curvature form of the normal bundle to $\nu_f$. Then by Chern-Weil
theory, we get bounds on Euler and Pontrjagin classes which
determine a vector bundle up to finitely many possibilities. An
alternative proof was recently found by K.~Tapp~\cite{Ta}. Instead
of getting bounds on the characteristic classes, Tapp estimates
the number of homotopy classes of maps into the classifying space.
\end{rmk}

\appendix
\section{Vector bundles with diffeomorphic total spaces}
The purpose of this appendix is to discuss to what extent a vector
bundle is determined by its total space. We got interested in this
problem when we noticed that under the assumptions of the
corollary~\ref{intro: no collapse on soul}, the homeomorphism
finiteness for the total spaces can be easily obtained from the
parametrized version of Perelman's stability theorem~\cite{Per1}.

Let $\eta_k$ be an infinite sequence of vector bundles over a
closed smooth manifold $M$ such that the total spaces $E(\eta_k)$
are homeomorphic. Assume that the natural homomorphism
$\mathrm{Homeo}(M)\to\mathcal E(M)$ of the homeomorphism group of
$M$ into the group of (free) homotopy classes of self-homotopy
equivalences of $M$ has finite cokernel. Then the homeomorphism
type of the total space determines (the topological equivalence
class of) a vector bundle, up to a finite ambiguity. (Indeed, the
homeomorphism $E(\eta_i)\to E(\eta_1)$ induces a self-homotopy
equivalence $g_i$ of $M$. Passing to a subsequence, we can assume
that $g_i^{-1}\circ g_j$ is homotopic to a homeomorphism for any
$i,j$. Let $s_i$ be the zero section of $\eta_i$. Then the maps
$s_j$ and $s_i\circ g_i^{-1}\circ g_j$ have equal invariants in
the sense of section~\ref{Invariants of maps}. Hence, for a fixed
$j$,~\ref{pe determine bundles} implies that the bundles
$(g_i^{-1}\circ g_j)^\#\eta_i$ fall into finitely many isomorphism
classes. Therefore, the bundles $\eta_k$ fall into finitely many
topological equivalence classes.)

For example, if the group $\mathcal E(M)$ is finite, then the
homeomorphism type of the total space determines (the isomorphism
class of) a vector bundle, up to a finite ambiguity. In fact,
since a vector bundle is determined, up to a finite ambiguity, by
its characteristic classes, it suffices to assume that the natural
action of $\mathcal E(M)$ on the cohomology groups that contain
these classes is an action of a finite group. For instance, up to
isomorphism, there are only finitely many vector bundles over
$CP^n$ with homeomorphic total spaces because $\mathcal E(CP_n)$
is finite
(see~\cite[ch12, 18.3]{Rut} for many more examples, also
see~\cite{DW} for the case when $M$ is a sphere).

Also if the image of the diffeomorphism group $\mathrm{Diffeo}(M)$
in $\mathcal E(M)$ has finite index, then
the homeomorphism type of the total space determines
(the smooth equivalence class of) a vector bundle, up to finite ambiguity.
If $dim(M)\ge 6$ and $\pi_1(M)=1$, then
the image of $\mathcal D(M)$ in $\mathcal E(M)$ is commensurable to
the isotropy subgroup of the total rational Pontrjagin class of
$TM$~\cite[pp322--323]{Sul}. In particular, the image of $\mathcal D(M)$ has
finite index in $\mathcal E(M)$ when $p(TM)=1$.
For instance, up to smooth equivalence, there are only finitely
many vector bundles with homeomorphic total spaces
over the direct product of finitely many spheres.

Now we give an example of an infinite
sequence of pairwise topologically inequivalent vector bundles
with diffeomorphic total spaces.
The base manifold will be homotopy equivalent to $S^3\times S^3\times S^5$.
We are thankful to Shmuel Weinberger for providing a key idea
for this example (as usual, the authors
assume all responsibility for possible mistakes).

\begin{ex}\label{append: ex}
Consider a manifold $N=S^3\times S^3\times S^5$ and let $t$ be a
nonzero element of $H^{8}(N,\mathbb Z)$. One can always find a
positive integer $q$ and a vector bundle $\tau$ over $N$ of rank
$\ge\dim(N)$ such that  $q t=p_2(\tau)\in H^{12}(N,\mathbb Z)$
where $p_2$ is the second integral Pontrjagin class. (Indeed, the
Pontrjagin character $ph$ defines an isomorphism of vector spaces
$$ph\co \widetilde{K}(N)\otimes\mathbb Q\to
\oplus_{i>0}H^{4i}(N,\mathbb Q).$$
Consider a finite sum $\sum_i
[\eta_i]\otimes p_i/q_i$, where $p_i$, $q_i$ are integers and
$\eta_i$ are bundles over $N$, such that $ph(\sum_i
[\eta_i]\otimes p_i/q_i)=t\otimes 1$.
Reducing  to a common
denominator and using that $p_i[\eta_i]$ is the class of
$p_i$-fold Whitney power of $\eta_i$, we get $t\otimes
1=ph(\xi)/q^\prime$ for some positive integer $q^\prime$ and a
bundle $\xi$ over $N$; one can choose the rank of $\xi$ to be any
number $\ge 15$. Since the first Pontrjagin class lives in the
zero group, the formula for the $ph_2$, the part of the Pontrjagin
character that lives in $8$th cohomology, reduces to $ph_2=-p_2/6$
and we are done.)

Replacing $\tau$ with its Whitney power, we can assume that
$\tau$ is fiber homotopy equivalent to the trivial bundle $TN$.
(Note that $p_2(\tau)$ is still an integral multiple of $t$.)
Since $dim(N)$ is odd and $\ge 5$, and $N$ is simply connected
a result of Browder-Novikov~\cite[II.3.1]{Br} implies that there is a
closed smooth manifold $M$ and a homotopy equivalence
$f\co M\to N$ such that  $f^\#\tau$ is stably isomorphic to $TM$.

It follows from~\cite{Sie} that any automorphism of
$H_3(S^3\times S^3,\mathbb Z)\cong\mathbb Z^2$
is induced by a self-homotopy equivalence.
Since the inclusion $S^3\times S^3\to N$ induces an isomorphism
of the third integral homology groups, any automorphism of
$H_3(N,\mathbb Z)$
is induced by a self-homotopy equivalence. The same is
therefore true for $M$.
Furthermore, any automorphism of
$H^{8}(M,\mathbb Z)$
is induced by a self-homotopy equivalence.
(Indeed, start with $\phi\in\mathrm{Aut}(H^{8}(M,\mathbb Z))$
and conjugate to by the Poincar\"e duality to get
an automorphism $\phi^\prime$ of $H_3(N,\mathbb Z)$.
If $f$ is a self-homotopy equivalence of $M$ inducing $\phi^\prime$,
then $\phi$ can be identified with the transfer map for $f$.
The transfer map is the inverse to $f^\ast$, hence
$f^{-1}$ induces $\phi$.)

Note that $\mathrm{Aut}(H^{8}(M,\mathbb Z))\cong GL(2,\mathbb Z)$.
Recall that any vector $(a,b)\in\mathbb Z^2$ is
$GL(2,\mathbb Z)$-equivalent to $(k,0)$ where $k=\gcd(a,b)$.
The vectors $v_m=(k,km)$ are $GL(2,\mathbb Z)$-equivalent to $(k,0)$
and lie in different orbits of the $GL(2,\mathbb Z)$-stabilizer of $(k,0)$.

Thus one can find an infinite sequence
of elements $w_m\in H^{8}(M,\mathbb Z)\cong \mathbb Z^2$
that lie in different orbits of
the stabilizer of $p_2(TM)$ in the group $GL(2,\mathbb Z)$
and are $GL(2,\mathbb Z)$-equivalent to $p_2(TM)$.
Find self-homotopy equivalences $g_m$ of $M$ such that
$w_m=g_m^\ast p_2(TM)$;
let $g_1=\mathrm{id}$.
Let $\eta_m$ be a bundle over $M$
with $[\eta_m]=[g_m^\#TM]-[TM]$; one can choose the rank of $\eta_m$
to be any number $\ge 15$.

Now if $\eta_i$ is topologically equivalent to
$\eta_j$, then there is a selfhomeomorphism $h$ of $M$
that takes $p_2(\eta_i)=p_2(g_i^\#TM)-p_2(TM)$ to
$p_2(\eta_j)=p_2(g_j^\#TM)-p_2(TM)$ ($p_2$ is additive because
$p_1=0$). Since homeomorphisms
preserve rational Pontrjagin classes,
we get $h^\ast g_i^\ast p_2(TM)=g_j^\ast p_2(TM)$, up
to elements of finite order.
The group $H^8(M,\mathbb Z)$ is torsion free, hence
the above equality hold exactly. So $h^\ast w_i=w_j$
which contradicts the definition of $w_m$. Thus, $\eta_m$ are
pairwise topologically inequivalent.

The map $g_m$ induces a homotopy equivalence
from the total space $E(\eta_m)$ to $E(\eta_1)$.
Now $[TE(\eta_m)]=[\eta_m]+[TM]=[g_m^\#TM]-[TM]+[TM]=[g_m^\#TM]$.
Thus, the homotopy equivalence $E(\eta_m)\to E(\eta_1)$
is tangential and hence is homotopic to a
diffeomorphism~\cite[pp226-228]{LS}.
\end{ex}

\begin{rmk}
We now briefly describe a generalization of the above
example. First of all, instead of $S^3\times S^3$,
one can start with a product $P$ of an arbitrary number of
$1$-connected homology $m$-spheres with odd $m>1$.
(Any odd-dimensional sphere is rationally homotopy equivalent
to a $K(\mathbb Z,m)$ space, so the product $P$ is rationally
homotopy equivalent $K(\mathbb Z^k,m)$ where $k$ is the number of factors.
This implies that the natural action
$\mathcal E(P)\to\mathrm{Aut}(H_m(P,\mathbb Z))\cong GL(k,\mathbb Z)$ has
finite cokernel. Thus, ``almost'' every automorphism is induced
by a homotopy equivalence which turns out to be enough for us.)

By taking a product with a suitable $1$-connected manifold $S$
(which was $S^5$ in the example) we can shift dimensions
so that $\dim(P\times S)$ is odd and
$$\mathcal E(P\times S)\to
\mathrm{Aut}(H^{4i}(P\times S,\mathbb Z))\cong GL(k,\mathbb Z)$$ has
finite cokernel for some $i$.
Using~\cite[I.3.1]{Br}, we replace $P\times S$ by a homotopy equivalent
manifold $M$ with a nonzero Pontrjagin class $p_i$.
(The freedom in the choice of $\tau$ gives infinitely many
possibilities for the homeomorphism type of $M$.)
Now it not hard to cook up an infinite sequence of
pairwise topologically inequivalent
bundles $\eta_j$ over $M$ with diffeomorphic total spaces.
One can also get nonsimply connected examples by taking
products of $E(\eta_j)\times L$ where $L$ is a suitable closed manifold.
\end{rmk}

\small
\bibliographystyle{amsalpha}
\bibliography{exp}

\end{document}